\documentclass{amsart}

\usepackage{a4wide}
\usepackage{amssymb}
\usepackage[utf8]{inputenc}
\usepackage[hyperref]{hyperref}

\renewcommand{\phi}{\varphi}
\renewcommand{\theta}{\vartheta}

\newcommand{\N}{{\mathbb{N}}}
\newcommand{\R}{{\mathbb{R}}}
\newcommand{\Q}{{\mathbb{Q}}}
\newcommand{\Z}{{\mathbb{Z}}}
\newcommand{\1}{{\mathbf{1}}}
\renewcommand{\epsilon}{\varepsilon}

\newcommand{\moduli}{{\mathcal{M}}}
\newcommand{\pairing}[2]{{\langle{#1}|{#2}\rangle}}
\newcommand{\abs}[1]{{\left\lvert #1\right\rvert}}
\newcommand{\norm}[1]{{\left\| #1\right\|}}
\newcommand{\POS}{{\mathbf{q}}}
\newcommand{\MOM}{{\mathbf{p}}}
\newcommand{\x}{{\mathbf{x}}}

\newcommand{\z}{{\mathbf{z}}}
\newcommand{\bg}{{\boldsymbol{\gamma}}}
\newcommand{\tu}{\tilde{u}}

\DeclareMathOperator{\CZ}{CZ}
\DeclareMathOperator{\evaluationmap}{ev}

\theoremstyle{plain}
\newtheorem{theorem}{Theorem}
\newtheorem{propo}[theorem]{Proposition}
\newtheorem{corollary}[theorem]{Corollary}

\theoremstyle{remark}
\newtheorem{remark}{Remark}

\theoremstyle{definition}
\newtheorem*{defi}{Definition}

\begin{document}

\title{Towards a good definition of algebraically overtwisted}

\author[F.~Bourgeois]{Fr\'ed\'eric Bourgeois}
\email[F.~Bourgeois]{fbourgeo@ulb.ac.be}
\address[F.~Bourgeois]{D\'epartement de Math\'ematique CP 218\\
  Universit\'e Libre de Bruxelles \\
  Boulevard du Triomphe\\1050 Bruxelles\\Belgium}

\author[K.~Niederkr\"uger]{Klaus Niederkr\"uger}
\email[K.~Niederkr\"uger]{kniederk@umpa.ens-lyon.fr}
\address[K.\ Niederkr\"uger]{\'Ecole Normale Sup\'erieure de Lyon\\
  Unit\'e de Math\'ematiques Pures et Appliqu\'ees\\
  UMR CNRS 5669\\ 46 all\'ee d'Italie\\ 69364 LYON Cedex 07\\
  France}

\begin{abstract}
  Symplectic field theory (SFT) is a collection of homology theories
  that provide invariants for contact manifolds.  We show that
  vanishing of any one of either contact homology, rational SFT or
  (full) SFT are equivalent.  We call a manifold for which these
  theories vanish \emph{algebraically overtwisted}.
\end{abstract}


\maketitle

\section{Introduction}

A (coorientable) contact structure $\xi$ on a $(2n+1)$--dimensional
manifold $M$ is a hyperplane field of the tangent bundle that can be
written as the kernel of a $1$--form $\alpha$ that satisfies the
inequality $\alpha\wedge d\alpha^n \ne 0$.  On closed manifolds
contact structures are stable under deformations, and their
equivalence classes are discrete sets.  Much effort has been invested
in understanding $3$--dimensional contact manifolds, and a rich theory
has been created.  For this, many different techniques have been
applied ranging from topological ones to different algebraic
invariants like Heegaard Floer theories \cite{HeegaardFloerContact}, 
or contact homology.  One of
the first basic properties that were discovered for $3$--manifolds was
the distinction between \emph{overtwisted} and \emph{tight} contact
structures.  Being overtwisted is a topological property, but it has
many consequences for algebraic invariants.

The algebraic invariant under consideration is Symplectic Field Theory (SFT),
introduced by Eliashberg, Givental and Hofer~\cite{SymplecticFieldTheory}. 
This large formalism 
contains in particular several versions of contact invariants, such as contact homology.
These invariants, described in Section \ref{sec:SFT}, are based on the count of
holomorphic curves in the symplectization of a contact manifolds and
are defined for contact manifolds of any odd dimension.

It has been proved by Eliashberg and Yau
\cite{YauContactHomologyVanishes} that the contact homology of any
overtwisted $3$--manifold is trivial.  Given that the classification of
such manifolds is purely topological \cite{Eliashberg_Overtwisted}, it
was to be expected that also the other invariants of SFT do not
provide interesting information.  In this article we confirm this
conjecture by a more general result for contact manifolds of any odd
dimension.  In fact, it follows already from purely algebraic
properties that vanishing of any of the different homology theories
implies that the other ones also have to be trivial.

\begin{theorem}\label{theorem algebraically overtwisted}
  Let $(M,\alpha)$ be a $(2n-1)$--dimensional closed contact manifold
  with a non degenerate contact form.  All of the following statements
  are equivalent:
  \begin{itemize}
  \item [(i)] The contact homology (without marked points) of
    $(M,\alpha)$ vanishes.
  \item [(ii)] The rational SFT (without marked points) of
    $(M,\alpha)$ vanishes.
  \item [(iii)] The SFT (without marked points) of $(M,\alpha)$
    vanishes.
  \item [(iv)] The contact homology with marked points of $(M,\alpha)$
    vanishes.
  \item [(v)] The rational SFT with marked points of $(M,\alpha)$
    vanishes.
  \item [(vi)] The SFT with marked points of $(M,\alpha)$ vanishes.
  \end{itemize}
\end{theorem}

\begin{remark}
  Note that any of these invariants may be defined over different
  coefficient rings.  In the theorem we assume that the same ring is
  used for the different homologies.
\end{remark}

To date no final generalization of overtwisted contact manifolds to
higher dimensions has been found.  This theorem, together with
\cite{YauContactHomologyVanishes}, motivates the following definition
that can be easily applied to any dimension.

\begin{defi}
  A contact manifold $(M,\alpha)$ is called \textbf{algebraically
    overtwisted}, if any of the homologies listed in
  Theorem~\ref{theorem algebraically overtwisted} vanishes.
\end{defi}

Several examples of algebraically overtwisted contact structures are
known: As stated above the contact homology of overtwisted manifolds
vanishes \cite{YauContactHomologyVanishes}.  Similarly Otto van Koert
and the first author of this article have extended this result by
showing that the contact homology of negatively stabilized contact
manifolds of any dimension is trivial
\cite{BourgeoisvKoertContactHomology}.

\begin{corollary}
  All of these examples are thus algebraically overtwisted, and have
  vanishing SFT.
\end{corollary}

A tentative generalization of the notion of overtwisted to higher
dimensions was given in \cite{NiederkruegerPlastikstufe}, where
$PS$--overtwisted manifolds were defined.  Current work by the
authors~\cite{PSOverwistedIsAlgebraicallyOvertwisted} will show that
$PS$--overtwisted manifolds also have vanishing contact homology.

\subsection*{Acknowledgments}
At the time of the creation of this article, we were both working at
the \emph{Universit\'e Libre de Bruxelles}.  The second author was
being funded by the \emph{Fonds National de la Recherche Scientifique}
(FNRS).

We thank Otto van Koert for fruitful discussions, and Hansj\"org
Geiges for valuable comments.

\section{Contact homology and variants of SFT}
\label{sec:SFT}

The different contact invariants above are homologies of certain
differential graded algebras with a $\1$--element, that means
\begin{defi}
  A differential graded algebra $(\mathcal{A},\partial)$ is a graded
  algebra, equipped with a differential $\partial:\,\mathcal{A}_* \to
  \mathcal{A}_{*-1}$ such that $\partial^2 = 0$ and which satisfies
  the graded Leibniz rule $\partial(a\cdot b) = (\partial a)\cdot b +
  (-1)^{\abs{a}} a\cdot (\partial b)$.
\end{defi}

The vanishing results in this article are all based on the following
easy remark.

\begin{remark}
  Let $(\mathcal{A},\partial)$ be a differential graded algebra
  with~$\1$, and denote the homology over that algebra by
  $H_*(\mathcal{A},\partial)$.  The homology vanishes if and only if
  $\1$ is an exact element, that means if there is an element $a\in
  \mathcal{A}$ such that $\partial a = \1$.
\end{remark}
\begin{proof}
  The $\1$--element is closed, because $\partial \1 = \partial
  (\1\cdot \1) = (\partial \1)\cdot \1 + (-1)^{\abs{\1}} \1 \cdot
  \partial \1= 2\,\partial \1$, hence it is obvious that $\1$ has to
  be exact for the homology to vanish.  On the other hand, if there is
  an element $a\in\mathcal{A}$ such that $\partial a = \1$, then the
  whole homology has to vanish, because we can write an arbitrary
  cycle $b\in\mathcal{A}$ as $b = \1\cdot b = (\partial a)\cdot b =
  \partial (a\cdot b) - (-1)^{\abs{a}} a\cdot \partial b = \partial
  (a\cdot b)$.
\end{proof}

Our proof of the main theorem is based on using algebraic properties
and exploiting that the contact homology algebra embeds naturally into
both the rational SFT algebra and the full SFT algebra.  In particular
the $\1$--elements all coincide under these inclusions.  Before
starting to describe the actual proof, we will briefly repeat how the
algebras are defined, and what the corresponding boundary operators
are (see \cite{SymplecticFieldTheory}).  Since the algebras of the
different versions of SFT are all build up by using closed Reeb
orbits, and the corresponding differentials all count certain
holomorphic curves, we will first fix some common notation.  Readers
familiar with symplectic field theory can safely skip the next
section, and skim back when needed.

\subsection*{Notation: Closed Reeb orbits and holomorphic curves}

Let $\alpha$ be a contact form for the $(2n-1)$--dimensional contact
manifold $(M,\xi)$.  The associated Reeb vector field $R_\alpha$ is
defined as the unique solution of the equations
\begin{equation*}
  \alpha(R_\alpha) = 1 \text{ and } i_{R_\alpha}d\alpha = 0 \;.
\end{equation*}
A closed Reeb orbit $\gamma$ of $R_\alpha$ is called non degenerate,
if the corresponding Poincar\'e return map does not have eigenvalues
of size $1$.  We call $\alpha$ a non degenerate contact form, if all
of its Reeb orbits are non degenerate.  Any contact form can be made
non degenerate by a small perturbation, and so we will always assume
from now on that $\alpha$ is non degenerate.  The associated Reeb
vector field $R_\alpha$ then has only countably many closed orbits,
and we can introduce a total order on the set of closed Reeb orbits
$\bg = (\gamma_1, \gamma_2,\dotsc)$.  Note that multiple orbits are
considered to be completely unrelated to the corresponding simple
orbits.  Denote the period of an orbit $\gamma$ by $T(\gamma)$, and
its multiplicity by $\kappa_\gamma$.  We fix a parametrization for
each closed orbit $\gamma_k$ by choosing a base point on $\gamma_k$.
A convenient short hand notation to handle ordered tuples of closed
Reeb orbits is to consider sequences $I = (i_k)_k\in\N_0^\N$ with only
finitely many non zero elements, and to denote by $\bg^I$ the tuple of
orbits
\begin{equation*}
  (\underbrace{\gamma_1, \dotsc, \gamma_1}_{i_1}, \dotsc,
  \underbrace{\gamma_N, \dotsc, \gamma_N}_{i_N}) \;,
\end{equation*}
where $N$ is large enough to capture all non vanishing elements of
$I$.  We allow for $I$ also the sequence $\textbf{0} = (0,\dotsc)$
giving rise to the empty tuple $\bg^\textbf{0} = ()$.  Finally, let
$\abs{I}$ be the number of non zero components in the sequence $I$,
and $C(I)$ be the integer
\begin{equation*}
  C(I) = \abs{I}!\,\, i_1!\dotsm i_N!\,
  \kappa_{\gamma_1}^{i_1}\dotsm \kappa_{\gamma_N}^{i_N}
\end{equation*}
again for $N$ large enough.

To compute the Conley-Zehnder index $\CZ(\gamma)$ of a closed Reeb
orbit $\gamma$, we have to fix a trivialization of the contact
structure $\xi$ along~$\gamma$.  To do this in a unified way, choose a
basis $A_1,\dotsc,A_s$ of $H_1(M,\Z)$ (for $H_1(M,\Z)$ with torsion,
we refer to \cite[Section~2.9.1]{SymplecticFieldTheory}) and for each
element $A_j$ a closed path $\psi_j$ representing $A_j$.  Fix a
trivialization of $\xi$ along $\psi_j$ in an arbitrary way, and choose
for every closed Reeb orbit $\gamma$, a surface $S_\gamma$ bounding
$\gamma$ and the corresponding combination of $\psi_j$'s that
represent $[\gamma]\in H_1(M,\Z)$, then use $S_\gamma$ to extend the
trivialization of $\xi$ from the $\psi_j$'s to $\gamma$.

Let $\bigl(\R \times M, d(e^t \alpha)\bigr)$ be the symplectization of
$(M,\alpha)$.  A complex structure $J$ on a symplectic vector bundle
$(E,\Omega)$ is called compatible with $\Omega$, if $\Omega(J\cdot,
J\cdot) = \Omega(\cdot, \cdot)$, and if $\Omega(J\cdot,\cdot)$ defines
a metric.  We choose a compatible $\R$--invariant complex structure
$J$ on the symplectic vector bundle $(\xi, d\alpha)$ and extend it to
an almost complex structure on the symplectization by $J
\frac{\partial}{\partial t} = R_\alpha$.  To define the differentials
for the different homologies, we have to enumerate certain holomorphic
curves.  Let $(\Sigma_g, j)$ be a compact Riemann surface of genus
$g$, and let $I^+ = (i_k^+)_k$, and $I^- = (i_l^-)_l$ be finite
sequences of integers.  Associate to every $i_k^+ \ne 0$ points
$\overline{x}^1_k,\dotsc,\overline{x}^{i^+_k}_k\in \Sigma_g$, to every
$i_l^-\ne 0$ points $\underline{x}_l^1,\dotsc,\underline{x}^{i_l^-}_l
\in \Sigma_g$, together with nonzero tangent vectors $\overline{v}_k^i
\in T_{\overline{x}_k^i} \Sigma_g$ and $\underline{v}_l^j \in
T_{\underline{x}_l^j} \Sigma_g$ respectively.  For reasons that will
become clear below, we call $\overline{\x} = \{ \overline{x}_k^l \}$
the positive, $\underline{\x} = \{ \underline{x}_k^l \}$ the negative
punctures, and the attached vectors are called asymptotic markers.
Additionally let there be $m$ marked points $y_1,\dotsc,y_m$ on the
Riemann surface $\Sigma_g$.  All the marked points and the positive
and negative punctures have to be pairwise distinct.

A map 
\begin{equation*}
  \tu=(a,u) :\, (\Sigma_g \setminus (\overline{\x} \cup \underline{\x}),
  j) \to (\R \times M, J)
\end{equation*}
is a $(j,J)$--holomorphic map, if $J\circ D\tu = D\tu\circ j$ for
every point of $\Sigma_g \setminus (\overline{\x} \cup
\underline{\x})$.  Additionally we require the following properties at
the punctures.  Choose for every puncture $p \in \overline{\x} \cup
\underline{\x}$ a holomorphic chart $D^2 \to \Sigma_g$ such the origin
is mapped to $p$, and such that the asymptotic marker points along the
positive real axis.  In polar coordinates $\bigl(\rho
e^{i\theta}\bigr)\in D^2$, the following asymptotic conditions have to
be satisfied by $\tu$:
\begin{equation*}
  \lim_{\rho \to 0} a\bigl(\rho e^{i\theta}\bigr) =
  \begin{cases}
    + \infty & \text{ if $p \in \overline{\x}$,} \\
    - \infty & \text{ if $p \in \underline{\x}$,}
  \end{cases}
  \quad \text{ and } \quad
  \lim_{\rho \to 0} u\bigl(\rho e^{i\theta}\bigr) =
  \begin{cases}
    \gamma_k(-\frac{T_k}{2\pi} \theta) & \text{ if  $p =
      \overline{x}_k^i$,} \\
    \gamma_l(\frac{T_l}{2\pi} \theta) & \text{ if $p =
      \underline{x}_l^j$,}
  \end{cases}
\end{equation*}
where $T_k$ denotes the period of the orbit $\gamma_k$.  When we do
not want to fix the complex structure on $\Sigma_g$, we call such a
map a $J$--holomorphic map.

Choose an additional puncture $x_0$ with asymptotic marker that will
be asymptotic to a closed Reeb orbit $\gamma$.  We denote by
$\moduli^A_{g, m} \bigl( \bg^{I^-};  \bg^{I^+}, \gamma\bigr)$ the space
of $J$--holomorphic maps as above that have an additional positive
puncture $x_0$, and by $\moduli^A_{g, m} \bigl(\gamma,
\bg^{I^-};\bg^{I^+}\bigr)$ the space of $J$--holomorphic maps that have an
additional negative puncture $x_0$.  In both cases, we assume that the
surface obtained by gluing $u(\Sigma_g \setminus (\{x_0\}\cup
\overline{\x} \cup \underline{\x}))$ with suitable surfaces
$S_{\gamma_k}$ represents the homology class $A \in H_2(M,\Z)$.

Let $(\Sigma_g,j)$ and $(\Sigma'_g,j')$ be compact Riemann surfaces
equipped with positive and negative punctures $\overline{\x},
\underline{\x}$ and $\overline{\x}', \underline{\x}'$ respectively and
with $m$ marked points $y_1,\dotsc,y_m$ and $y'_1,\dotsc,y'_m$.  We
call a diffeomorphism $\phi:\Sigma_g\to\Sigma'_g$ a reparametrization,
if it is a biholomorphism that is compatible with all special points.
This means that $\phi$ satisfies the equation $\phi_*j = j'$, and
$\phi(y_k) = y'_k$, and corresponding relations for the ordered
punctures.  The map also has to respect the asymptotic markers at each
puncture.  Define an equivalence relation $\sim$ on the space of maps
$\moduli^A_{g, m} (\dotsc)$ by saying that two maps $\tu=(a,u)$ and
$\tu'=(a',u')$ are equivalent, if there is a shift $\tau\in\R$,
and a reparametrization $\phi:\, (\Sigma_g,j) \to (\Sigma_g,j')$ such
that
\begin{equation*}
  (a, u) = (a'\circ\phi + \tau,u'\circ \phi)\;.
\end{equation*}

The moduli spaces
\begin{equation*}
  \widehat \moduli^A_{g, m} (\dotsc) =
  \moduli^A_{g, m} (\dotsc)/\sim
\end{equation*}
are obtained by dividing out the corresponding space of maps by the
equivalence relation $\tu\sim \tu'$ just defined.  Denote the first Chern
class of the complex vector bundle $(\xi = \ker\alpha,J)$ by
$c_1(\xi)$.  If the elements of $\widehat \moduli^A_{g, m} (\dotsc)$
are not branched coverings, then choosing $J$ generically these moduli
spaces are smooth orbifolds of dimension
\begin{align*}
  \begin{split}
    \dim \widehat\moduli^A_{g, m} \bigl(\bg^{I^-}; \bg^{I^+}, \gamma
    \bigr) &=
    (n-3)\,\left(2 - 2g -\abs{I^-} -\abs{I^+}-1\right) - 1 + 2m + \CZ(\gamma)\\
    & + 2\, \pairing{c_1(\xi)}{A} + \sum_{j=1}^\infty (i^+_j -
    i^-_j)\,\CZ(\gamma_j)
  \end{split} \\
  \begin{split}
    \dim \widehat\moduli^A_{g, m} \bigl(\gamma,
    \bg^{I^-}; \bg^{I^+} \bigr) &=
    (n-3)\,\left(2 - 2g -1 -\abs{I^-} -\abs{I^+}\right) - 1 + 2m  - \CZ(\gamma)\\
    & + 2\, \pairing{c_1(\xi)}{A} + \sum_{j=1}^\infty (i^+_j -
    i^-_j)\,\CZ(\gamma_j)
  \end{split} 
\end{align*}
that are equipped with a smooth evaluation map at the marked points
\begin{equation*}
  \evaluationmap:\, \widehat\moduli^A_{g, m}
  (\dotsc) \to M^{m},\,
  [a,u] \mapsto \bigl(u(y_1),\dotsc, u(y_m)\bigr) \;.
\end{equation*}
The moduli spaces have compactifications ${\overline \moduli}^A_{g, m}
\bigl(\bg^{I^-}; \bg^{I^+},\gamma\bigr)$, and ${\overline
  \moduli}^A_{g, m} \bigl( \gamma, \bg^{I^-};\bg^{I^+}\bigr)$
respectively consisting of holomorphic buildings of arbitrary height
\cite{BourgeoisCompactness}.

In the presence of branched coverings, the new ongoing approach to
transversality by Cieliebak and Mohnke (see \cite{CM} for the
symplectic case) or the polyfold theory developed by Hofer, Wysocki
and Zehnder \cite{H,HWZ} give to the moduli space $\overline
\moduli^A_{g, m} (\dotsc)$ the structure of a branched manifold (with
rational weights) with boundary and corners.  The presence of these
rational weights is due to the use of multivalued perturbations.

In the absence of marked points, and when $\dim \widehat \moduli^A_{g,
  0} (\dotsc) =0$, this moduli space consists of finitely many
elements with rational weights.  We denote the sum of these rational
weights by $n^A_g\bigl(I^-; I^+,\gamma\bigr)$ or
$n^A_g\bigl(\gamma,I^-; I^+\bigr)$.  When $m \ne 0$, we define a
multilinear form $n^A_{g,m}(\dotsc)$ on $m$--tuples of closed
differential forms $\Theta_1, \dotsc, \Theta_m$ on $M$ by the formula
\begin{equation*}
  \pairing{n^A_{g,m}(\dotsc)}{(\Theta_1, \dotsc, \Theta_m)} = 
  \int_{\overline \moduli^A_{g, m} (\dotsc)}
  \evaluationmap^*(\Theta_1 \times \dotsm \times \Theta_m) \;.
\end{equation*}
By convention, we set the multilinear form to $0$, if $\sum
\deg\Theta_j \ne \dim \moduli^A_{g, m} (\dotsc)$, and we define a
$0$--multivalued form $n^A_{g,0}(\dotsc)$ just by using the sum
$n^A_g(\dotsc)$ of rational weights defined above.

To define the algebras, we have to find a suitable coefficient ring.
For this choose a submodule
\begin{equation*}
  \mathcal{R} \le \bigl\{A \in H_2(M,\Z) \bigm|\,
  \pairing{c_1(\xi)}{A} = 0\bigr\}
\end{equation*}
to construct the group ring $\Q\bigl[H_2(M,\Z)/\mathcal{R}\bigr]$,
whose elements will be written as $\sum_{j=1}^k c_j e^{A_j}$, where
$c_j\in\Q$ and $A_j\in H_2(M,\Z)/\mathcal{R}$.  Different choices of
$\mathcal{R}$ may lead to different SFT invariants.  We define a
grading on $\Q\bigl[H_2(M,\Z)/\mathcal{R}\bigr]$ by $\abs{c\, e^{A}} =
-2\,\pairing{c_1(\xi)}{A}$.

Associate to every closed Reeb orbit $\gamma$ the formal variables
$q_\gamma$ and $p_\gamma$ with gradings
\begin{equation*}
  \abs{q_\gamma} = \CZ(\gamma) + n - 3 
  \quad \text{ and }
  \abs{p_\gamma} = -\CZ(\gamma) + n - 3
\end{equation*}
Given a finite sequence of integers $I$, we denote by $\POS^I$ the
monomial $q_{\gamma_1}^{i_1}\dotsm q_{\gamma_N}^{i_N}$ and by $\MOM^I$
the monomial $p_{\gamma_1}^{i_1}\dotsm p_{\gamma_N}^{i_N}$ for $N$
large enough.

Consider a formal variable $\hbar$ with grading $\abs{\hbar} =
2\,(n-3)$.

\subsection*{Contact homology}

The algebra $\mathcal{A}_{CH}$ of contact homology consists of
polynomials in the $q_\gamma$'s with coefficients in
$\Q\bigl[H_2(M,\Z)/\mathcal{R}\bigr]$.  Every element can be written
as a finite sum
\begin{equation*}
  f = \sum_{k=1}^K f_k \;,
\end{equation*}
where each term $f_k$ is of the form
\begin{equation*}
  f_k = c_k\,e^{A_k} \POS^{I_k}
\end{equation*}
with every $c_k\,e^{A_k} \in \Q\bigl[H_2(M,\Z)/\mathcal{R}\bigr]$ and
$I_k = \bigl(i_{j,k}\bigr)_j$ is a sequence of the type described
above.  The grading for each such monomial is given by
\begin{equation*}
  \abs{f_k} = \sum_{j=1}^\infty \Bigl(\CZ(\gamma_j) 
  + n-3\Bigr)\,i_{j,k} - 2\,\pairing{c_1(\xi)}{A_k} \;.
\end{equation*}

The sum of monomials $f = c_1 e^{A}q_{\gamma_1}^{i_1}\dotsm
q_{\gamma_N}^{i_N}$, and $g = c_2 e^{B}q_{\gamma_1}^{j_1}\dotsm
q_{\gamma_N}^{j_N}$ (we assume $N$ to be large enough to include all
non zero terms of both sequences $I = (i_k)_k$ and $J = (j_k)_k$) is
formal, and the multiplication of $f$ and $g$ gives
\begin{equation*}
  fg =  c_1c_2 e^{A+B}q_{\gamma_1}^{i_1}\dotsm
  q_{\gamma_N}^{i_N}q_{\gamma_1}^{j_1}\dotsm
  q_{\gamma_N}^{j_N} \;,
\end{equation*}
where we still have to permute the $q$--variables to get a monomial in
normal form. For this, we impose supercommutativity $q_{\gamma}
q_{\gamma^\prime} = (-1)^{|q_\gamma|\,|q_{\gamma^\prime}|}
q_{\gamma^\prime}q_{\gamma}$.  The differential on this algebra
$\mathcal{A}_{CH}$ is defined by
\begin{equation*}
  \partial q_{\gamma} = \sum_{A,I} \frac{n^A_{0,0}\bigl(
    \bg^I;\gamma\bigr)}{C(I)}\,e^A \,\POS^I \;,
\end{equation*}
where the sum runs over all integer valued sequences $I$ and all
homology classes $A$.  Remember that $n^A_{0,0}\bigl(
\bg^I;\gamma\bigr)$ counts (in the sense defined above) punctured
holomorphic spheres with a single positive puncture asymptotic to
$\gamma$, and negative punctures asymptotic to the orbits in $\bg^I$.
Effectively, the sum in the definition of the differential operator is
finite.  On one hand, the period of $\gamma$ gives an upper bound for
$\sum i_k T(\gamma_k)$ (see for example
\cite[Lemma~5.16]{BourgeoisCompactness}), so that only finitely many
sequences $I$ need to be taken into account.  On the other hand, the
compactness theorem for the space $\cup_A \moduli^A_{g, m} (
\bg^I;\gamma)$ with fixed $I$ shows that there may only be holomorphic
curves for finitely many choices of $A$.

For products extend the differential according to the graded Leibniz
rule, i.e. $\partial (fg) = \left(\partial f\right)\, g +
(-1)^{\abs{f}}\,f\, \left(\partial g\right)$.

\subsection*{Rational SFT}

The algebra $\mathcal{A}_{rSFT}$ of rational SFT can be interpreted as
a Poisson algebra with a distinguished element $\mathbf{h}$. Since we
are just interested in showing that its homology vanishes, we will
only describe it as a differential graded algebra.  The elements of
$\mathcal{A}_{rSFT}$ can be written as
\begin{equation*}
  f = \sum_{I^+}  f_{I^+}(\mathbf{q})\, \MOM^{I^+} \;,
\end{equation*}
where the sum runs over all finite sequences $I^+$ of integers, and
the coefficients $f_{I^+}(\mathbf{q}) \in \mathcal{A}_{CH}$ are
elements in the contact homology algebra that depend on $I^+$.  In
other words, the elements of $\mathcal{A}_{rSFT}$ are formal power
series in $p$--variables with coefficients in the contact homology
algebra $\mathcal{A}_{CH}$.  The grading of a monomial is given by
\begin{equation*}
  \abs{c\,e^{A}
    \POS^{I^-} \MOM^{I^+}} = \abs{c\,e^{A}\POS^{I^-}} 
  + \sum_{j=1}^\infty \bigl(n-3- \CZ(\gamma_j)\bigr)\,i^+_j \;.
\end{equation*}
The product between variables is supercommutative
\begin{equation*}
  q_\gamma q_{\gamma'} = (-1)^{\abs{q_\gamma} |q_{\gamma'}|}
  q_{\gamma'} q_\gamma,\quad
  q_\gamma p_{\gamma'} = (-1)^{\abs{q_\gamma}|p_{\gamma'}|}
  p_{\gamma'} q_\gamma \quad
  \text{ and }\quad p_\gamma p_{\gamma'} =
  (-1)^{\abs{p_\gamma}|p_{\gamma'}|}
  p_{\gamma'} p_\gamma \;.
\end{equation*}
The differential $d^\mathbf{h}$ is defined on a single $q$--variable
by the formula
\begin{align*}
  d^\mathbf{h} q_\gamma &= \sum_{A,I^-,I^+} \frac{n_{0,0}^A
    \bigl(\bg^{I^-};\bg^{I^+},\gamma\bigr)}{C(I^-)\,C(I^+)}\, e^A
  \POS^{I^-}\MOM^{I^+} \intertext{and on a $p$--variable respectively
    by} d^\mathbf{h} p_\gamma &= (-1)^{|p_\gamma|+1} \sum_{A,I^-,I^+}
  \frac{n^A_{0,0}\bigl( \gamma,\bg^{I^-};\bg^{I^+}
   \bigr)}{C(I^-)\,C(I^+)}\,e^A \POS^{I^-}\MOM^{I^+} \;.
\end{align*}
We are summing over all combinations of monomials $e^A
\POS^{I^-}\MOM^{I^+}$.  For the definition of the rational numbers
$n_{0,0}^A(\dotsc)$, we refer to Section ``Notation: Closed Reeb
orbits and holomorphic curves''.  For arbitrary elements in
$\mathcal{A}_{rSFT}$ extend the operator $d^\mathbf{h}$ by using the
graded Leibniz rule.

\subsection*{Full SFT}

The algebra of symplectic field theory $\mathcal{A}_{SFT}$ is composed
of formal power series of the form
\begin{equation*}
  F =  \sum_{g=0}^\infty \sum_{I^+}  f_{g,I^+}(\POS)\,
  \MOM^{I^+}\,\hbar^g \;,
\end{equation*}
as above $f_{g,I^+}(\POS)$ is an element in the contact homology
algebra, and $\hbar$ is a new formal variable of degree $2\,(n-3)$, so
that the total degree of a monomial $c
e^{A}\,\POS^{I^-}\MOM^{I^+}\hbar^g$ is given by
\begin{equation*}
  \abs{c\,e^{A}\,\POS^{I^-}\MOM^{I^+}\hbar^g} =
  \abs{e^{A}\,\POS^{I^-}\MOM^{I^+}} + 2g\,(n-3) \;.
\end{equation*}
Unlike $\mathcal{A}_{CH}$ and $\mathcal{A}_{rSFT}$, the algebra
$\mathcal{A}_{SFT}$ is not supercommutative, but instead has the
following commutator relations.  For two $q$--variables or two
$p$--variables the commutator relations are identical to the ones of
rational SFT, but for mixed terms, we require the supersymmetric
relation
\begin{equation*}
  [q_\gamma,p_{\gamma'}] :=
  q_\gamma p_{\gamma'} - (-1)^{\abs{q_\gamma}|p_{\gamma'}|} p_{\gamma'}
  q_\gamma =
  \begin{cases}
    \kappa_\gamma\hbar & \text{ if $\gamma = \gamma'$} \\
    0 & \text{ otherwise.}
  \end{cases}
\end{equation*}
One way to incorporate this relation into a formalism is by
representing $p_\gamma$ as the derivation operator
$\kappa_\gamma\hbar\, \partial/\partial q_\gamma$.  The
$\hbar$--variable commutes with the $q$-- and $p$--variables.

In \cite{SymplecticFieldTheory}, the differential $d^\mathbf{H}$ of
SFT was given as the commutator with a distinguished element
$\mathbf{H}$, but here we will just specify the effect of
$d^\mathbf{H}$ on the generators of $\mathcal{A}_{SFT}$
\begin{align*}
  d^\mathbf{H} \hbar & = 0\;, \\
  d^\mathbf{H} q_\gamma &= \sum_{g, A, I^-, I^+} \frac{n^A_{g,0}
    \bigl(\bg^{I^-}; \bg^{I^+},\gamma\bigr)}{C(I^-)\,C(I^+)}\, e^A
  \POS^{I^-}\MOM^{I^+}\hbar^g\;, \intertext{and} d^\mathbf{H} p_\gamma
  &= (-1)^{|p_\gamma|+1} \sum_{g, A, I^-, I^+} \frac{n^A_{g,0} \bigl(\gamma,\bg^{I^-}; 
    \bg^{I^+}\bigr)}{C(I^-)\,C(I^+)}\, e^A \POS^{I^-}\MOM^{I^+}\hbar^g
  \;,
\end{align*}
and extend it to general elements by the graded Leibniz rule.

\subsection*{Marked points}

Choose closed differential forms $\Theta^1,\dotsc,\Theta^d$ that
represent an integral basis for the de Rham cohomology ring
$H^*_\mathrm{dR}(M)$.  Take any of the algebras described above, i.e.,
let $\mathcal{A}$ be either the contact homology algebra
$\mathcal{A}_{CH}$, the rational symplectic field algebra
$\mathcal{A}_{rSFT}$, or the SFT algebra $\mathcal{A}_{SFT}$.  Define
new formal variables $t_1,\dotsc,t_d$ with grading $\abs{t_j} :=
\deg\Theta^j - 2$.  Let $\mathcal{A}^*$ be the algebra of formal power
series in the $t_1,\dotsc,t_d$ with coefficients in $\mathcal{A}$ such
that the $t_j$ supercommute among themselves, and with the
$q$--variables, and possibly also (if they are part of $\mathcal{A}$)
with the $p$--, and $\hbar$--variables.

The differential $\partial^*$ on $\mathcal{A}^*$ vanishes on $\hbar$,
and the $t_1,\dotsc,t_d$
\begin{equation*}
  \partial^* t_j  = 0,\text{ and } \partial^* \hbar  = 0\;.
\end{equation*}

To define $\partial^*$ on the $q$-- and $p$--variables, introduce
first the notation $\boldsymbol\Theta = \sum_{j=1}^d \Theta^jt_j$, and
\begin{equation*}
  \pairing{n^A_{g,m} (\dotsc)}{(\boldsymbol\Theta, \dotsc,
    \boldsymbol\Theta)} = \sum_{1\le a_1, \dotsc, a_m\le d}
  \pairing{n^A_{g,m} (\dotsc)}{(\Theta^{a_1}, \dotsc,
    \Theta^{a_m})}\,t_{a_1}\dotsm t_{a_m}\;,
\end{equation*}
and then set
\begin{align*}
  \partial^* q_\gamma &= \sum_{g,m, A, I^-, I^+} \frac{\pairing{n^A_{g,m}
      \bigl(\bg^{I^-}; \bg^{I^+},\gamma\bigr)}{(\boldsymbol\Theta,
      \dotsc, \boldsymbol\Theta)}} {C(I^-)\,C(I^+)}\, e^A
  \POS^{I^-}\MOM^{I^+} \,\hbar^g \;, \intertext{and}
  \partial^* p_\gamma &= (-1)^{|p_\gamma|+1} \sum_{g,m, A, I^-, I^+} \frac{\pairing{n^A_{g,m}
      \bigl(\gamma, \bg^{I^-}; 
      \bg^{I^+}\bigr)}{(\boldsymbol\Theta, \dotsc,
  \boldsymbol\Theta)}}
  {C(I^-)\,C(I^+)}\, e^A \POS^{I^-}\MOM^{I^+}\, \hbar^g \;,
\end{align*}
where we are summing over all $I^-$ described above, and in case
$\mathcal{A} = \mathcal{A}_{CH}$, we assume that $g=0$, and $I^+ =
\mathbf{0}$.  If $\mathcal{A} = \mathcal{A}_{rSFT}$, we still keep
$g=0$, but allow any sequence $I^+$, and finally if $\mathcal{A} =
\mathcal{A}_{SFT}$, any integer $g\ge 0$, and sequence $I^+$ is
allowed.

\section{Proof of Theorem~\ref{theorem algebraically overtwisted}}

\subsection{The implications $(ii) \Rightarrow (i)$, $(iii)
  \Rightarrow (i)$, and $(n+iii) \Rightarrow (n)$}

The statements are based on the following trivial remark.

\begin{remark}
  Let $\pi:\,\mathcal{A}' \to \mathcal{A}$ be a chain map between two
  differential graded algebras $(\mathcal{A},\partial)$ and
  $(\mathcal{A}',\partial')$.  From $\pi\circ\partial' = \partial
  \circ \pi$ it follows immediately that an exact element $f'\in
  \mathcal{A}'$ is mapped to an exact element $f = \pi(f') \in
  \mathcal{A}$.
\end{remark}

Corresponding to each of the cases $(ii) \Rightarrow (i)$, $(iii)
\Rightarrow (i)$, and $(n+iii) \Rightarrow (n)$, we find an element
$g$ in either $\mathcal{A}_{rSFT}$, $\mathcal{A}_{SFT}$, or
$\mathcal{A}^*$ such that $d^\mathbf{h} g = \1$, $d^\mathbf{H} g =
\1$, or $\partial^* g = \1$.

For the first case, define a projection $\pi:\, \mathcal{A}_{rSFT} \to
\mathcal{A}_{CH}$ by mapping any monomial
\begin{equation*}
  f_{I^+}(\POS)\, \MOM^{I^+} \mapsto
  \begin{cases}
    f_{I^+}(\POS) & \text{ if $I^+ = \mathbf{0}$} \\
    0 & \text{ otherwise,}
  \end{cases}
\end{equation*}
and extending this map linearly.  It is clear that $\pi$ will be an
algebra homomorphism, and to see that it is a chain map, just compare
the definitions of the differentials $\partial$ and $d^{\mathbf{h}}$.
All terms counted by the contact homology differential also appear in
the rational SFT differential, and it is clear that $\partial\circ \pi
= \pi\circ d^{\mathbf{h}}$ holds, if $d^{\mathbf{h}}$ cannot decrease
the number of $p$--variables in any monomial.  Note that already by
the Leibniz rule, the differential $d^\mathbf{h}$ can decrease the
number of factors in a monomial at most by one, and to decrease the
number of $p$--factors, variables $p_\gamma$ have to exist such that
$d^\mathbf{h}p_\gamma$ contains (non zero) terms without any
$p$--coordinates at all.  This would only be possible, if there were
non empty moduli spaces $\moduli^A_{0,0}\bigl(\bg^\mathbf{0}; \gamma,
\bg^{I^-}\bigr)$ of spheres without positive punctures, but by the
maximum principle no such curves exist.  We show in Appendix~\ref{sec:
  maxprin} that a weaker form of the maximum principle still holds for
solutions of perturbed Cauchy-Riemann equations, so that the same
conclusion remains true.  It follows that $\pi \bigl(d^\mathbf{h}
f_{I^+}(\POS)\, \MOM^{I^+}\bigr) = 0$ for any monomial with $I^+\ne
\mathbf{0}$, and so $\pi(g)\in\mathcal{A}_{CH}$ will be a primitive of
the unit element $\1$.

The chain map $\pi:\, \mathcal{A}_{SFT} \to \mathcal{A}_{CH}$ for the
second case will be defined similarly, dropping any monomial that
contains a positive $\hbar$-- or $p$--power.  This map is compatible
with the commutator relations, and it is again a chain map, because
$d^\mathbf{H}$ cannot decrease the number of $p$--factors without
raising the $\hbar$--power and vice versa such that $\pi
\bigl(d^\mathbf{H} f_{I^+}(\POS)\, \MOM^{I^+}\hbar^g\bigr) = 0$, if
either $g\ne 0$ or $I^+\ne \mathbf{0}$, and so the contact homology
algebra is trivial as we wanted to show.

To prove $(n+iii)\Rightarrow (n)$, let $\pi:\,\mathcal{A}^* \to
\mathcal{A}$ be the projection that drops any monomial containing a
$t_j$--variable.  Let $\partial$ be the differential of $\mathcal{A}$
(i.e., depending on $\mathcal{A}$ either $\partial$, $d^\mathbf{h}$ or
$d^\mathbf{H}$).  We need to show $\pi$ is a chain map.  As before the
argument here is that $\partial$ is the zero order term of
$\partial^*$ in the $t_j$, and this is true because the count $n^A_{g}
(\dotsc)$ coincides by definition with
$\pairing{n^A_{g,0}(\dotsc)}{()}$, furthermore $\partial^*$ can never
decrease the number of $t_j$--variables of a monomial, so that
$\partial \pi f = \1$.

\subsection{The implications $(i) \Rightarrow (ii)$, $(i) \Rightarrow
  (iii)$, and $(n) \Rightarrow (n+iii)$}

For the implication $(i) \Rightarrow (ii)$, assume that $f_0$ is an
element in the contact homology algebra $\mathcal{A}_{CH}$ such that
$\partial f_0 = \1$.  This $f_0$ canonically embeds into the rational
SFT algebra $\mathcal{A}_{rSFT}$, where we can compute its
differential $d^\mathbf{h} f_0 = \1 - g$.  The element
$g\in\mathcal{A}_{rSFT}$ is always closed, because $d^\mathbf{h} g =
d^\mathbf{h}(\1- d^\mathbf{h}f) = 0$, and all of its terms contain at
least one $p$--variable.  The formal inverse of $\1 - g$ is given by
\begin{equation*}
  (\1-g)^{-1} := \sum_{k=0}^\infty g^k\;,
\end{equation*}
where we set $g^0 = \1$.  This object is well defined in
$\mathcal{A}_{rSFT}$, because all terms in $g$ have at least one
$p$--factor, so that only the powers of $g$ up to $k=\abs{I^+}$
contribute to the terms $\MOM^{I^+}$ in $(1-g)^{-1}$, and in
particular $(\1-g)^{-1}$ is a formal power series in the
$p$--variables.  It is obvious that $d^\mathbf{h}(1-g)^{-1} = 0$.

Define an element $f\in\mathcal{A}_{rSFT}$ by
\begin{equation*}
  f := f_0\,(\1-g)^{-1}\;.
\end{equation*}
It easily follows that
\begin{equation*}
  d^\mathbf{h} f = d^\mathbf{h}\left(f_0\,(\1-g)^{-1}\right)
  = (d^\mathbf{h}f_0)\,(\1-g)^{-1} - f_0\,d^\mathbf{h}(1-g)^{-1}
  = (\1-g)\,(\1-g)^{-1} = \1 \;,
\end{equation*}
just as we wanted to show.

We will now prove the implication $(i)\Rightarrow (iii)$ in a similar
way: Assume that $f\in\mathcal{A}_{CH}$ is such that $\partial f =
\1$.  We can canonically embed the contact homology algebra
$\mathcal{A}_{CH}$ into the algebra of symplectic field theory
$\mathcal{A}_{SFT}$, since the commutation relations for only
$q$--variables are identical in both spaces.  As we said above,
dropping any term with a $p$-- or $\hbar$--variable the differential
$d^\mathbf{H}$ coincides on $\mathcal{A}_{CH}$ with the differential
$\partial$ of contact homology.  Thus
\begin{equation*}
  d^\mathbf{H} f = \partial f - G = \1 -G \;,
\end{equation*}
where $G$ only contains monomials with a non zero power of $\hbar$ or
of $p$, therefore the formal inverse $(\1 -G)^{-1} = \sum_{k=0}^\infty
G^k$ is a well defined element in $\mathcal{A}_{SFT}$.  Moreover it is
closed, and so
\begin{equation*}
  F = f\, \sum_{k=0}^\infty G^k
\end{equation*}
is a primitive of $\1$, because by using the Leibniz rule
\begin{equation*}
  d^\mathbf{H} \left(f\,\sum_{k=0}^\infty G^k\right) =
  \left(d^\mathbf{H} f\right) \,\sum_{k=0}^\infty G^k
  - f\, \sum_{k=0}^\infty d^\mathbf{H} \bigl(G^k\bigr) = \1\;,
\end{equation*}
follows, and the homology of SFT vanishes.

Finally, we compare the invariants with and without marked points, and
prove $(n)\Rightarrow (n+iii)$: Let thus $\mathcal{A}$ be either
$\mathcal{A}_{CH}$, $\mathcal{A}_{rSFT}$, or $\mathcal{A}_{SFT}$
without marked points, and $\mathcal{A}^*$ the corresponding algebra
with marked points.  Use that $\partial$ and $\partial^*$ are
identical in $0$--th order of $t_j$--powers, so that if $\partial f =
\1$, we have that $\partial^* f = \1 - G$, where all terms in $G$ have
positive $t_j$--powers.  As above, the formal inverse $(1-G)^{-1} =
\sum G^k$ is an element of $\mathcal{A}^*$, so that we can define $F
:= f\,(1-G)^{-1}$ which is a primitive of $\1$ with respect to
$\partial^*$.

\appendix

\section{Maximum principle for perturbed holomorphic curves}
\label{sec: maxprin}

Let $T_0 > 0$ be smaller than the period of any closed Reeb orbit in
$M$.  Let $\widehat G$ be the positive (but not definite) metric on
$\R \times M$ defined by $\widehat G(\cdot, \cdot) = d\alpha(\cdot, J
\cdot)$.  Let $G$ be the positive definite metric on $\R \times M$
such that the Reeb field $R_\alpha$ and the Liouville field $
\frac{\partial}{\partial t}$ are mutually orthogonal, are orthogonal
to $\xi$, and have unit length.  In the next proposition, $\nu \in
\Lambda^{0,1}(\Sigma_g, \tu^*T(\R \times M))$ will denote a perturbation
for the Cauchy-Riemann equation.

\begin{propo}
  A curve $\tu:\, \Sigma_g \setminus ( \x \cup \z ) \to \R \times M$
  that satisfies the perturbed Cauchy-Riemann equation
  \begin{equation*}
    d\tu + J \circ d\tu \circ j = \nu
  \end{equation*}
  with $\norm{\nu}_{L^2(G)} < 2 \sqrt{T_0}$ has to have top punctures
  $\x \neq \emptyset$.
\end{propo}

\begin{proof}
  Assume that there exists a map $\tu :\, \Sigma_g \setminus \z \to \R
  \times M$ satisfying $d\tu + J \circ d\tu \circ j = \nu$ that is
  asymptotic for $t \to -\infty$ to the orbits $\gamma_1, \dotsc,
  \gamma_s$.  By Stokes theorem, we have
  \begin{equation*}
    \int_{\Sigma_g \setminus \z} \tu^* d\alpha =
    - \sum_{i=1}^s T_i < - T_0 \;,
  \end{equation*}
  where $T_i$ is the period of $\gamma_i$.
 
  On the other hand, let $z = x+iy$ be coordinates of a complex chart
  on $\Sigma_g$ such that $\norm{\partial_x} = \norm{\partial_y} = 1$.
  Then
 \begin{equation*}
   d\alpha\bigl(\partial_x \tu, \partial_y \tu\bigr) =
   -d\alpha\bigl(J \partial_y \tu, \partial_y \tu\bigr)
   + d\alpha\bigl(\nu(\partial_x), \partial_y \tu\bigr)
   = \widehat G\bigl(\partial_y \tu, \partial_y \tu\bigr) +
   \widehat G\bigl(\nu(\partial_x), -J\,\partial_y \tu\bigr) \;.
 \end{equation*}
 After a suitable rotation in the $(x,y)$--plane, we obtain, using
 Cauchy-Schwarz inequality, and $\bigl(\norm{d\tu}_{\widehat G} -
 1/2\,\norm{\nu}_{\widehat G}\bigr)^2 \ge 0$
 \begin{equation*}
   d\alpha\bigl(\partial_x \tu, \partial_y \tu\bigr)
   = \norm{d\tu}_{\widehat G}^2
   + \widehat G\bigl(\nu(\partial_x), -J\,\partial_y \tu\bigr)
   \ge  \norm{d\tu}_{\widehat G}^2 - \norm{\nu}_{\widehat G}\, \norm{d\tu}_{\widehat G}
   \ge -\frac{1}{4}\, \norm{\nu}_{\widehat G}^2 \;,
 \end{equation*}
 where $\norm{\cdot}_{\widehat G}$ is the semi-norm induced by
 $\widehat G$.

Integrating over $\Sigma_g \setminus \z$, we obtain
\begin{equation*} 
  \int_{\Sigma_g \setminus \z} \tu^* d\alpha \ge
  -\frac{1}{4}\, \norm{\nu}_{L^2(\widehat G)}^2 \;.  
\end{equation*}
Comparing with the Stokes bound for this integral, we obtain
$\norm{\nu}_{L^2(G)} \ge \norm{\nu}_{L^2(\widehat G)} > 2 \sqrt{T_0}$,
a contradiction.
\end{proof}


\end{document}